\documentclass[12pt]{amsart}
\usepackage{epsfig}

\usepackage{graphicx}
\usepackage{amscd}
\usepackage{amsmath}
\usepackage{amsfonts}
\usepackage{amssymb}
\usepackage{verbatim}

\newtheorem{theorem}{Theorem}[section]
\theoremstyle{plain}

\newtheorem{corollary}[theorem]{Corollary}

\newtheorem{defi}[theorem]{Definition}

\newtheorem{prop}[theorem]{Proposition}

\def\bdot{\mbox{\bf\,.\,}}

\def\wtil{\widetilde}
\def\til{\tilde}
\def\wt{\wtil}
\def\what{\widehat}
\def\wQ{\what{Q}}
\newcommand{\clos}{\mbox{\rm clos}}

\newcommand{\lam}{\lambda}

\newcommand{\gam}{\gamma}

\newcommand{\sig}{\sigma}

\newcommand{\R}{{\mathbb R}}

\newcommand{\Z}{{\mathbb Z}}

\def\N{{\mathbb N}}

\def\ba{{\bf a}}

\def\Max{{\sf Max}}
\def\Min{{\sf Min}}

\def\ov{\overline}

\begin{document}

\thispagestyle{empty}

\title[Extension of the Morse transformation]{The adic realization of the Morse transformation and the extension
of its action on the solenoid}

\author{A. M. Vershik}

\author{B. Solomyak}

\address{Anatoly Vershik, Mathematical Institute of Russian Academy
of Sciences, St.Petersburg branch, Fontanka 27, St.Petersburg,
191023, Russia}
\email{vershik@pdmi.ras.ru}
\address{Boris Solomyak, Box 354350, Department of Mathematics,
University of Washington, Seattle WA 98195}
\email{solomyak@math.washington.edu}

\date{\today}

\begin{abstract}
We consider the adic realization of the Morse transformation on the
additive group of integer dyadic numbers. We discuss the arithmetic
properties of that action. Then we extend that action to the action
of the group of rational dyadic numbers on the solenoid.
\end{abstract}
\maketitle

\rightline{To the memory of Alexander Livshits}

Sasha Livshits (1950-2008) was the author of one of the most important
 theorems of modern dynamics, which is well-known now, --- the theorem about
cohomology of hyperbolic systems. He proved it when he was a student.
Later he worked on many other problems of symbolic dynamics,
ergodic theory and combinatorics. His deep and important ideas made a great impression on those who
interacted  with him (this includes the second author).
 The first author considers him the best of his students.

\section{Introduction}

The Morse dynamical system was discovered by Morse and popularized
by Hedlund-Gottshalk. Later is was studied by many authors (see
\cite{Quef,Pyth} and references there), as a simplest nontrivial
substitution. Moreover, it was historically the first example of a
substitution. It is generated by the Thue-Morse sequence, which was
extensively studied from the point of view of combinatorics of words
and symbolic complexity (see \cite{AlSh,Moshe} and references
there). The new approach to symbolic dynamics and ergodic
transformations (adic transformation) which was suggested by the
first author \cite{VerDan}, can be applied also to substitutions
(stationary adic transformations); this was done in the paper of A.
Livshits and the first author (\cite{VerLiv}). Later, other authors
developed this connection in the context of topological dynamics,
see \cite{Forrest,Dur}, but here our focus is on measure-preserving
transformations. The adic realization of a substitution dynamical system allows one
to consider simultaneously not only the substitution itself but also
the one-sided shift which accompanies any substitution. The idea
which was advocated  by the first author in \cite{Ver.tran} is to
consider the natural extension of that shift and correspondingly
extend the substitution system in order to make an
 essential link between the theory of substitutions and hyperbolic dynamics.
In this paper we consider the ``two-sided extension'' of the Morse system
 which yields {\it the Morse action of the countable group $Q_2$}
 (the group of dyadic rational numbers)
 on the group of its characters --- solenoid $\wQ_2$, reworking more carefully
 and correcting some details of \cite{Ver.tran}.
 We obtain also some new properties of the adic realization of the Morse transformation.
 One of the corollaries of the adic approach is an explicit calculation which shows
 how to obtain the Morse system as a time-change of the dyadic odometer.
 The operation of differentiation of dyadic sequences plays an
 important role in our constructions.
 The spectral theory of the Morse system, which goes back to Kakutani
 \cite{Kaku} (see also \cite{Quef,Pyth} and references there)
 is also becoming more transparent under these considerations, but we
 do not address it in this paper.

In Section 2 we collect a series of auxiliary results, which use the
adic realization of the Morse system; some of them new, others are
well-known, but we give  their adic version. We discuss in more
detail (than in  \cite{Ver.tran}) the so-called ``Morse arithmetic".
Section 3 describes the two-sided extension of the Morse
transformation and embedding of it to the Morse action of the group
$Q_2$ on the group $\wQ_2$.  We also formulate some problems.

One should consider this article as an attempt, looking at the
special case of the Morse transformation, to attack the general
problem of defining a two-sided extension of a substitution system,
and corresponding to this extension its embedding into an action of
a larger group. The final goal of the constructions is to show the
link between the theory of substitutions and of hyperbolic systems.

The first author thanks to Max Plank Institute for Mathematical
Physics in Bonn and Erwin Schrodinger International Institute for
Mathematical Physics in Vienna where this article was finished, and
Russian funds RFBR-08-01-000379a, NSh-2460.20081112.
The research of the second author was supported in part by NSF grant DMS-0654408.

\section{Definitions and  adic realizations of Morse system.}

\subsection{Morse as a substitution system.} Consider the alphabet $\{0,1\}$.
The Morse substitution is defined by $\zeta(0)=01,\ \zeta(1)=10$; it
is extended to words in $\{0,1\}$ by concatenation. The Thue-Morse
sequence (sometimes also called Prouhet-Thue-Morse sequence) is a
fixed point of the substitution
\begin{equation} \label{def-Morse}
u = u_0u_1u_2\ldots=\lim_{n\to\infty} \zeta^n(0)=0110100110010110\ldots
\end{equation}
It has many remarkable features,
see e.g.\ \cite{AlSh} and \cite[Ch.2,\,Ch.5]{Pyth}. It is easy to see that
$$
u[0,2^{n+1}-1] = u[0,2^n-1]\,\ov{u[0,2^n-1]}\ \ \mbox{for}\ n\ge 0,
$$
where we denote $u[i,j]=u_i\ldots u_j$ and
$\ov{w}$ is the ``flip'' of a word $w$ in the alphabet $\{0,1\}$, that is, we interchange
$0 \leftrightarrow 1$. The sequence
$u$ is non-periodic, but uniformly recurrent, with well-defined uniform
frequencies of subwords. It is also known
that $u_n$ is the sum of the digits (mod 2) in the binary representation of $n$.

Let $\sigma$ be the left shift on the spaces $\{0,1\}^\N$ and
$\{0,1\}^\Z$ with the product topology.
 The substitution dynamical system is sometimes
considered on the space of one-sided sequences, and sometimes on the
space of two-sided sequences. \footnote{We have used above the
notion ``one-sided" and ``two-sided" in completely different sense,
see also below.}

The ``one-sided'' substitution space is defined as the orbit closure
of $u$ under the shift:
$$
X^+_\zeta = \clos\{\sig^n u:\,n\ge 0\}.
$$
The ``two-sided'' substitution space is defined as the set $X_\zeta$
of all bi-infinite sequences in $\{0,1\}^\Z$ whose every block
(subword) occurs in $u$. The substitution dynamical systems are
$(X^+_\zeta,\sig)$ (one-sided) and $(X_\zeta,\sig)$ (two-sided). The
advantage of the two-sided system is that it is a homeomorphism,
whereas the one-sided is not. Measure-theoretically these two
systems are isomorphic: in fact, both are minimal and uniquely
ergodic, and the one-sided system is a.e.\ invertible. The
restriction of any point in $X_\zeta$ to non-negative coordinates is
in $X^+_\zeta$, and all but countable many elements of $X^+_\zeta$
have a unique extension to points in $X_\zeta$. The exceptions are
$u$, which extends to $u\bdot u$ and $\ov{u}\bdot u$, and its
``flip'' $\ov{u}$, as well as their orbits.

\subsection{Adic realization.} Here we follow the general definition of the
adic transformation from \cite{VerDan} and \cite{VerLiv}, but focus
only on the Morse system as it was done in \cite{Ver.tran}. Consider
$\textbf{Z}_2\cong\{0,1\}^\N$, the compact additive group of 2-adic
integers, and the odometer (``adding machine'') transformation $T$,
which is adic transformation by definition - this is group
translation on $\textbf{Z}_2$ (see below).  We obtain the adic
realization of Morse transformation by {\it changing the order of
symbols $0,1$ depending on the next symbol}. Namely, consider the
lexicographic order on $\textbf{Z}_2$ induced by
$$
0 \prec_0 1,\ 1 \prec_1 0
$$
as follows:
$$
\{x_i\} \prec \{y_i\} \ \ \Longleftrightarrow\ \ \exists\,j:\ x_i=y_i \ \mbox{for} \ i>j\ \mbox{and}\
x_j \prec_z y_j, \ \mbox{where} \ z=x_{j+1}=y_{j+1}.
$$
This is a partial order; two sequences are comparable if they are
cofinal (agree except in finitely many places). The set of maximal
points is $\Max = \{(01)^\infty,(10)^\infty\}$, and the set of
minimal points is $\Min = \{(0)^\infty, (1)^\infty\}$. \footnote{We
denote infinite periodic sequence with period $(ab\ldots c)$ as
$(ab\ldots c)^{\infty}$.} Let $M$ be the immediate successor
transformation in the order $\prec$ on $\textbf{Z}_2$. Here we write
down the formulas for the action of $M$ explicitly. If
$x\not\in\Max$, then $x$ starts with $(01)^n00$ or $(01)^n1$, or
$(10)^n0$, or $(10)^n11$, where $n\ge 0$. We have
\begin{equation} \label{eq-M1}
\begin{array}{lr}
M((01)^n00*) =(1^{2n+1} 0*), &  M((01)^n1*) =(0^{2n}1*),\\
M((10)^n0*) =(1^{2n}0*), & M((10)^n11*) =(0^{2n+1} 1*) \end{array}
\end{equation}
Note that it is well-defined everywhere except on the two maximal
points: $\Max = \{(01)^\infty,(10)^\infty\}$. It is easy to see that
$M$ is continuous on $\textbf{Z}_2 \setminus \Max$. But it is
impossible to extend $M$ by continuity to those points: there are no
well-defined limits of $\lim_{n\to \infty} M((01)^n*)$, and
$\lim_{n\to \infty} M((10)^n*)$, because
$$\lim_{n\to \infty} M((01)^n00*)=(1)^{\infty},$$
but
$$\lim_{n\to \infty} M((01)^n1*)=(0)^{\infty}.$$
Analogously,
$$\lim_{n\to \infty} M((10)^n00*)=(1)^{\infty},$$
but
$$\lim_{n\to \infty} M((10)^n1*)=(0)^{\infty}.$$
Since we also have two minimal points, we can extend $M$ to a
bijection arbitrarily, setting
\begin{equation} \label{eq-Minf}
M((01)^\infty) = (1)^\infty,\ \ \ M((10)^\infty) = (0)^\infty,
\end{equation}
or vice versa. But this extension is not continuous  at those
points.

The obvious corollary of the definition of $M$ is that it commutes
with ``flips,'' that is,
\begin{equation} \label{eq-flip}
M(\ov{x}) = \ov{M(x)}\ \ \mbox{for all}\ x\in \textbf{Z}_2.
\end{equation}
The action of $M$ on $\textbf{Z}_2 \setminus \Max$ may
be expressed as follows: we scan the sequence $x$ to the right
until we see two identical symbols $aa$ and replace
the beginning of the sequence by $\ov{a}\ldots \ov{a} a$, keeping the second
occurrence of $a$ and everything that follows unchanged.

\subsection{The relation of adic model to traditional representation}
Now we indicate the relation between the dynamical systems
$(\textbf{Z}_2,M)$ and $(X_\zeta,\sig)$. Let
$$
g:\ \textbf{Z}_2 \to X_\zeta,\ \ \ \ g(x) = \{(M^{n-1} x)_0\}_{n\in \Z}.
$$
We have the following diagram:
$$
         \begin{CD}
         \textbf{Z}_2     @>M>> \textbf{Z}_2 \\
         @VVg V       @VVg V\\
         X_\zeta     @>\sig >> X_\zeta
       \end{CD}
$$
It is obvious from the definition that the diagram commutes.
It is also easy to see that $g$ is surjective,
continuous on $\textbf{Z}_2\setminus \wt{\Max}$, and $g(0^\infty) = \ov{u}\bdot u$.
Here we denoted by $\wt{\Max}$ the set of points in ${\bf Z}_2$
which are cofinal with the points in $\Max$ (equivalently, their left semi-orbit).

It may be useful to write down $g^{-1}$ explicitly. Consider the substitution map on $X_\zeta$:
$$
\zeta:\ X_\zeta\to X_\zeta,\ \ \zeta(\ldots a_{-2}a_{-1}\bdot a_0 a_1\ldots) =
\ldots \zeta(a_{-2})\zeta(a_{-1})\bdot \zeta(a_0)\zeta(a_1)\ldots
$$
It is well-known (and easy to see) that for every $\ba\in X_\zeta$
there is a unique $\ba'\in X_\zeta$ such that either $\ba = \zeta(\ba')$ or $\ba = \sig\zeta(\ba')$, and these
cases are mutually exclusive.
Let $\Psi:\ X_\zeta \to X_\zeta$ be given by $\Psi(\ba) = \ba'$. Then we have the following commutative diagram:
$$
         \begin{CD}
        \textbf{Z}_2     @>\sig >> \textbf{Z}_2 \\
         @VVg V       @VVg V\\
         X_\zeta     @>\Psi >> X_\zeta
       \end{CD}
$$
where $\sig$ is the left shift on $\textbf{Z}_2$. Therefore, to compute the $n$-th symbol of $g^{-1}(\ba)$ we need to
take $(\Psi^n(\ba))_0$ for $n= 0,1,2\ldots$

Now let us explain why this model is richer than the ``two-sided"
model $X_\zeta$. In the adic realization, we have the adic
transformation, which is isomorphic, up to neglecting of two orbits,
to the substitution, AND we have the one-sided shift in the space
$\textbf{Z}_2$. Evolution under the adic transformation of the first
digit $x_0$ of the sequence $\{x_i\}\in \textbf{Z}_2$ just gives the
orbit of $u$ under the transformation $\sig$ on $ X_\zeta$. The
one-sided shift in the space ${\bf Z}_2$ in terms of substitutions
is the substitution itself --- i.e. the transformation of any
sequence which substitutes $0$ by $01$ and $1$ by $10$. So we have
in the adic model a simultaneous realization of both
transformations: the adic and the shift (substitution). The problem
arises how to include in this picture the natural extension of the
one-sided shift
--- the two-sided shift, and at the same time how to extend the adic transformation
to the whole space. We will do this in the next section but first we interpret a familiar property
of the Morse system in our terms.

\subsection{Morse as a 2-point extension of the odometer.}
 \begin{defi}The classical 2-odometer is the following affine transformation $T$ on
 the additive group ${\bf Z}_2$ of dyadic integers:
 $$T:\, Tx=x+1.$$
 \end{defi}
 The transformation $T$ preserves the Haar (=Bernoulli, Lebesgue)
 measure on the group $\textbf{Z}_2$.
It is well-known that the Morse system can be represented as a group
(2-point) extension of the dyadic odometer. This is the most popular
point of view on the Morse system in dynamics. The adic realization
of the Morse transformation gives another way to look at this fact;
it is an illustration of the importance of the so-called
differentiation of sequences.

 Define the important map:
  \begin{defi}Differentiation of sequences is the map $D:\,{\bf Z}_2\to {\bf Z}_2$:
$$
  D(\{x_n\}_{n=0}^{\infty}) = \{(x_{n+1}-x_n)\ \mbox{\rm mod}\ 2,\ n=0,1, \dots\}.
$$
\end{defi}
This is nothing other than a 2-to-1 factorization of $\textbf{Z}_2$
on itself. It is clear that the differentiation commutes with the
``flip" which has been defined above: $D(\ov{x}) = D(x)$.

In spite of simplicity of the definition of the map $D:\,{\bf Z}_2
\to {\bf Z}_2$, there are no good and simple ``arithmetic" or
``analytic" expressions for the description of $D(\cdot)$. Recently V.
Arnold for different reasons made many experiments on the behavior
of 0-1 sequences under the iteration of differentiation \cite{A}.
But the most important thing for us is that the map $D$ takes the
Morse transformation into the odometer:
\begin{prop}
The following equality takes place: $T\circ D = D\circ M$.
\end{prop}
This is an immediate corollary of (\ref{eq-M1}). So in the adic
realization, the Morse transformation $M$ is a 2-covering of the
odometer in its algebraic form. Let us give a precise description of
the equivalence between the Morse transformation and a 2-extension
of the odometer. Let $F(x) = (Dx, x_0)$ be the map from
$\textbf{Z}_2$ to $\textbf{Z}_2 \times \{0,1\}$. This is a
bijection, and we have the following commutative diagram:
$$
         \begin{CD}
       \textbf{Z}_2     @>M >> \textbf{Z}_2 \\
         @VV F V       @VV F V\\
       \textbf{Z}_2\times \{0,1\}    @>T(\phi) >>  \textbf{Z}_2\times \{0,1\}
       \end{CD}
$$
Here $T(\phi)$ is a 2-extension of $T$ with a cocycle $\phi$ on $Z_2$ defined by
\begin{equation} \label{def-cocycle}
\phi(y)=\left\{\begin{array}{ll} 0 & \mbox{if $y$ starts with odd number of 1's}\\
                                 1 & \mbox{if $y$ starts with even number of 1's}
\end{array} \right.
\end{equation}
To make it work on maximal elements we also need to set $\phi(1^\infty) = \phi(0^\infty)=1$.
Recall that the group extension is defined by
$$
T(\phi)(x,g)  = (Tx, \phi(x) + g).
$$
We have
$$
M=F^{-1}T(\phi)F,
$$
so $M$ is canonically isomorphic to a 2-extension of the odometer
$T$ with cocycle $\phi$. We can identify $\textbf{Z}_2$ with
$\textbf{Z}_2 \times \{0,1\}$ regarding the second component
$\{0,1\}$ as a new digit of a sequence; then the map $F$ becomes a
transformation of $\textbf{Z}_2$, and we can consider the extension
inside the group $\textbf{Z}_2$.
We give another interpretation of this cocycle in the next section.

\medskip

{\bf Remark.} The Morse system can also be realized as a 2-point extension of the odometer
in the traditional substitution form, and it is interesting that the projection is again
given by the differentiation map. This follows from the fact that for the Thue-Morse sequence $u$
(see (\ref{def-Morse})), its derivative sequence $D(u) = 1011101010\ldots$ is the
fixed point of the substitution $0\to 11,\ 1 \to 10$ (see \cite[p.201]{AlSh}),
which generates a measure-preserving transformation isomorphic to the 2-odometer.

\medskip

Denote by $S$ the map $x \rightarrow 2x$ on  $\textbf{Z}_2$, this is
nothing other than the one-sided non-invertible shift or Bernoulli
endomorphism if we represent the elements of
$\textbf{Z}_2$ as a sequence of $0$'s and $1$'s. It is easy to check the
following fact.
\begin{prop}
The 2-odometer as well as the Morse transformation satisfy the
following equation:
$$TS=ST^2, \quad  MS=SM^2.$$
\end{prop}

\subsection{Morse as a time change of the odometer and Morse arithmetic.}
Because the group of rational integers $\Z$ is a dense invariant subgroup of the group of dyadic
integers, we can consider the Morse transformation
$M$ in the adic realization as a map of integers to itself.
This subsection is based on \cite[p.538]{Ver.tran}, but we provide more details.

Let us identify a sequence $x_0x_1 x_2\ldots$ with the dyadic decomposition
of a number: $\sum_j x_j 2^j$. Here is the list of several first values of $M(n)$:

\begin{center} \begin{tabular}{ccccccccccccccccc}
0 & 1 & 2 & 3 & 4 & 5 & 6 & 7 & 8 & 9 & 10 & 11 & 12 & 13 & 14 & 15 & \ldots \\ \hline
1 & 3 & 7 & 2 & 5 & 15 & 4 & 6 & 9 & 11 & 31 & 10 & 13 & 8 & 12 & 14 & \ldots \end{tabular}
\end{center}
\medskip
The table can be easily verified using (\ref{eq-M1}).

To this end, we introduce the following sequence:
$$
a_r = \left\{ \begin{array}{cc} \frac{2^r-1}{3}\, & \mbox{if}\ r\equiv 0\ (\mbox{mod}\ 2), \\[1.4ex]
                                \frac{2^r-2}{3}\, & \mbox{if}\ r\equiv 1\ (\mbox{mod}\ 2).\end{array} \right.
$$
Each $n\in \N$ is uniquely represented in one of the following ways ($r=r(n)$):
\begin{equation} \label{eq-ord1}
n= \left\{ \begin{array}{ll}  2^r \ell + a_{r-1} & \mbox{(i)}, \\
                              2^r \ell + 2^{r-1} + a_r & \mbox{(ii)}, \end{array} \right.
\end{equation}
where $\ell\ge 0$ is an integer.
  Define the mapping $M:\,\N \to \N\setminus \{0\}$ by
\begin{equation} \label{eq-ord2}
M(n)= \left\{ \begin{array}{ll}  n+a_{r(n)} & \mbox{in the case (i)}, \\
                                 n-a_{r(n)} & \mbox{in the case (ii)}. \end{array} \right.
\end{equation}
Although these formulas look a bit mysterious, they easily follow from (\ref{eq-M1}). In fact,
$$
a_r=\frac{2^r-1}{3} = (10)^{(r-2)/2},\ \ r\equiv 0\ (\mbox{mod}\ 2), \ \ \ \
a_r=\frac{2^r-2}{3} = (01)^{(r-1)/2},\ \ r\equiv 1\ (\mbox{mod}\ 2).
$$
Case (i) above occurs when the first pair $aa$ in the binary representation of $n$ is $00$.
Then $M(n)$ replaces the beginning of the sequence with $1$'s, which increases the number by $a_{r(n)}$
(observe that $a_{r-1}+a_r = 2^{r-1}-1=(1)^{r-1}$ independent of the parity of $r$). Similarly,
the case (ii) above occurs when the first pair  $aa$  in the binary representation of $n$ is
$11$. In this case $M(n)$ decreases or increases the number $n$ by $a_{r(n)}$.

Thus, we have described independently the restriction of the adic Morse system to $\N$;
that is why we use the same symbol $M$.
Define $M$ for negative integers by $M(-n) = -M(n-1)-1$. Then it is easy to check that we have the properties
$M(\ov{x}) = \ov{M(x)}$ where $\ov{n}=-n-1$; this should be understood by identifying integers with their
binary expansions. Thus, we have $M:\,\Z\to \Z \setminus \{0,-1\}$.
Note that $0 = (0)^\infty$ and $-1 = (1)^\infty$ are
the two minimal points in our ordering on ${\bf Z}_2$. According to (\ref{eq-Minf}),
$$
M(-1/3) \equiv M((10)^{\infty}) = (0){^\infty}\equiv 0\ \ \mbox{and}\ \ \ M(-2/3) \equiv M((01)^{\infty})
= (1)^{\infty}\equiv -1.
$$
\subsection{The orbit equivalence of Morse system and 2-odometer}
The orbit of the point $x \in X$ with respect to an invertible
transformation $S$ of $X$ is a set of $S^n x,\ n\in \Z$. Evidently,
the $T$-orbit of any point $x\in {\bf Z}_2$ which has infinitely many
$0$'s and $1$'s is the set of all points which are eventually equal
to $x$. The set of all points which have finitely many $0$'s or
$1$'s make one orbit (this is the common $T$-orbit of
$(0)^{\infty}$ and $(1)^{\infty}$). Let us describe the orbit
partition of Morse transformation, which follows directly from the
definition (2).
\begin{prop}
If a point $x\in {\bf Z}_2$ has infinitely many subwords $00$ and
subwords $11$, then the $M$-orbit of $x$ is the set of all points
which are eventually equal to $x$; the remaining countable set of the
points which have finitely many subwords $00$ or subwords $11$, is
just the union of four  semiorbits of $M$: two positive $M$-semiorbits
(${\Z}_+$ $M$-orbit)--- of the point $(0)^\infty$ and point
$(1)^\infty$, and two negative $M$-semiorbits (${\Z}_-$ $M$-orbit)
--- of the point $(10)^{\infty}$ and point $(01)^{\infty}.$
\end{prop}
Note that the negative $M$-semiorbit of $(10)^{\infty}$
(correspondingly $(01)^{\infty}$) consists of the set of points
which are eventually $(10)^{\infty}$ (correspondingly
$(01)^{\infty}$) and have an initial {\it even} word.

\begin{corollary}Orbit partitions of the 2-odometer and the adic realization of
the Morse transformation coincide  $(\!\!\mod 0)$ with respect to the Haar (Lebesgue)
measure on ${\bf Z}_2$.
\end{corollary}
As we saw, they coincide on the complement of a countable set. We
will refine this claim below.

Using our extension  of $M$ by the definition (3) we can make an additional
remark about those four semiorbits; we do not use it later. Note
that two positive $M$-semiorbits generate one $T$-orbit, and each
negative $M$-semiorbit is a full $T$-orbit. So in our definition (3)
we cut one $T$-orbit of $(0)^{\infty}$ and $(1)^{\infty}$ and glue
the first part $(0^{\infty})$ with the $M$-semiorbit of the point
$(10)^{\infty}$, and the second part $(1^{\infty})$ with the
$M$-semiorbit of the point $(01)^{\infty}$.

\medskip

 If $x \in \N \subset {\bf Z}_2 $ then we can
write tautologically
$$M(n)=T^{M(n)-n}(n),$$ where in the left-hand side $M(n)$ is the
image of $n$ under the transformation $M$, and in the right-hand
side $M(n)$ is a natural number. Now observe that by the definition
of the action of the Morse automorphism $M$ on the set of integers
defined by (\ref{eq-ord2}) above we have: $$M(n)-n=(-1)^{r(n)} \cdot
{a_{r(n)}}.$$

It is worth mentioning here that the value of the cocycle $\phi(n)$ from the previous subsection is exactly
$M(n)-n$ (mod 2), e.g. it is equal to $0$ iff $n$ and $M(n)$ have the same parity.

Denote $$\theta(n)=(-1)^{r(n)}\cdot a_{r(n)},$$ then we have the formula
$$M(n)=T^{\theta(n)}n$$ for each rational integer $n$. It is clear that the
function $r(\cdot)$ and consequently the function $\theta(\cdot)$ can be extended from positive integers $\N$
to the group of all dyadic integers ${\bf Z}_2$
as follows: formula (5) makes sense for all $x \in {\bf Z}_2$ with some $r\in \N$ and $\ell
\in {\bf Z}_2$, not only for integers $x$,
with the same definition. We just consider infinite sequences of $x_n$.
So $\theta(\cdot)$ becomes a function on ${\bf Z}_2$ with integer values; we can say that this
is simply the extension of $\theta(\cdot)$ by continuity in the pro-2-topology.

We have proved the following:

\begin{theorem}
Let $M$ be the adic realization of the Morse transformation in the space ${\bf Z}_2$.
Let $\wt{\Max}\cup \wt{\Min}$ be the countable set which is union of the semi-orbits under the action of $M$
of the  four points of ${\bf Z}_2$ : $$(0)^\infty,\ \  (1)^\infty,\ \ (01)^\infty,\ \ (10)^\infty.$$
Then on the $M$-invariant set ${\bf Z}_2\diagdown(\wt{\Max}\cup \wt{\Min})$ the odometer
$T:\ Tx=x+1$ and the Morse transformation $M$
have the same orbit partition, and moreover,
$$Mx=T^{\theta(x)}x \quad\mbox{for}\quad x\in {\bf Z}_2\diagdown (\wt{\Max}\cup \wt{\Min}),$$
where $\theta(x)$ is the function defined above.
\end{theorem}

The formula above gives an independent definition of the Morse transformation using a time change of the odometer.

Dye's Theorem asserts that any ergodic automorphism $S$ is isomorphic (mod $0$) to an
automorphism which is a time change of the odometer $T$ (or any other given ergodic automorphism):
$Sx=T^{\theta(x)}(x)$.
Nevertheless, there are few examples of an explicit formula for such a time change
function $\theta(\cdot)$.
The theorem above is just of this type: the Morse automorphism is represented as a time change
of the dyadic odometer. It is also known (see \cite{B}, theorem 3.8) that the time change
integer-valued function $\theta(\cdot)$
 cannot be integrable if ergodic automorphisms have the same orbits, unless $T =  S$ or  $T = S^{-1}$.
It is easy to check that our function $\theta$ is indeed non-integrable because it has exactly
two singularities on the space ${\bf Z}_2$ at the points $(01)^{\infty}(\equiv-1/3)$, and $(10)^{\infty}(\equiv-2/3)$,
and the measure of the cylinder on which the values of $\theta(x)$ are equal to $a_r$ is of order $C2^{-r}$, hence the
singularities have the type of simple poles $1/t$. The weakness (closeness to integrability)
of these singularities shows that the Morse automorphism is in a sense very close to the odometer,
i.e.\ to an automorphism with discrete spectrum.

\medskip

\textbf{Question.} What is the group generated by two
transformations of $\textbf{Z}_2$ --- the odometer $T$ and the Morse
transformation $M$? Is it a free group?

\section{Extension of the Morse transformation up to action of the group $Q_2$
on the solenoid}

In this section we define the so-called two-sided extension of the Morse transformation which
acts on the group of characters of dyadic rational numbers. It is an elaboration of
\cite[p.539]{Ver.tran} with important changes and additions.

\subsection{Preliminary facts about dyadic groups $Q_2$, $\what{Q}_2$, etc.}

Consider the exact sequence
$$
1\longrightarrow \Z \longrightarrow Q_2 \longrightarrow Q_2/\Z \longrightarrow 1
$$
where $Q_2$ is the countable additive group of real dyadic rational
numbers ($r/2^m, r\in \Z, m\geq 0$), the subgroup $\Z \subset Q_2$
is the group of rational integers, and the quotient group $Q_2/\Z$
is the group of all the roots of unity of orders $2^n,\ n=0,1, \dots
$ (a subgroup of the rotation group).

The group $Q_2$ can be presented as an inductive limit
$$\lim_{{\longrightarrow}_n} (\Z, w_n),$$ of the groups $\Z$,
with the embedding of $n$-th group given by
$$w_n(x)=2x,\quad  n=0,1 \dots.$$

  Consider the corresponding dual exact sequence for the groups of characters of the groups above:
  $$
 1\longleftarrow \R/\Z \longleftarrow \what{Q}_2 \longleftarrow \textbf{Z}_2\longleftarrow 1
 $$
  The group of characters of the group $Q_2/\Z$ is just the additive group of dyadic integers,
  $\textbf{Z}_2$,  which we considered in the previous sections, and which is the  inverse limit of
$2^n$-cyclic groups:
    $$\textbf{Z}_2=\lim_{\longleftarrow}(\Z/{2^n},p_n),$$
   with the maps $p_n:\Z/{2^n}\rightarrow \Z/2^{n-1},\ p_n(x)=x \ \mbox {mod}\ 2^{n-1}$.
 The group of characters of the group $\Z$ is the rotation group $S^1=\R/\Z$ (or the unit circle).

 Our main object --- the group $\what{Q}_2$ of characters of the group $Q_2$, is the so-called
 {\it 2-solenoid} and can be presented as an inverse limit of the rotation groups:
 $$\what{Q}_2=\lim_{{\longleftarrow}_n}(\R/\Z, v_{n+1}),\  n=0,1, \dots $$ where
  the homomorphisms are $$v_n: \R/\Z \rightarrow \R/\Z,\ \ v_n (u)=2u,\ n=1,2,\dots.$$

The group $\textbf{Z}_2$ is a closed subgroup of the group
$\what{Q}_2$ of those elements which have the zero projection on
$\R/\Z$ equal to $0$.

  The additive group $Q_2$ of dyadic rational
numbers is naturally embedded into $\what{Q}_2$ as a dense subgroup;
it consists of those characters which send the elements of the group
$Q_2$ to the roots of unity of degree $2^n$.

 Note that the additive group of the locally compact field $\textbf{Q}_2$ of all 2-adic number
 is naturally embedded into the solenoid  $\what{Q}_2$ as a dense subgroup of those elements of
$\what{Q}_2$, which have the projection under the map
$\what{Q}_2\rightarrow  \R/\Z$ to a root of unity
 of degree $2^n$ for some $n$:
 $$\textbf{Q}_2 \subset \what{Q}_2.$$

 As a compact group, $\what{Q}_2$ has a normalized Haar measure which is the product of
 Haar measures on the groups $\textbf{Z}_2$ and $\R/\Z$.

 The group {$\what{Q}_2$  is not the direct product of the subgroups $\textbf{Z}_2$ and
 the quotient group $\R/\Z$; there is a nontrivial 2-cocycle $\textbf{Z}_2$ of the group $\R/\Z$
 with values in $\Z$ (integer part of product).
But there is a non-algebraic decomposition into
 the direct product of the subgroup and the quotient group:
 $$
 \what{Q}_2 =\R/\Z \times \textbf{Z}_2.
 $$
  Here the first component of the decomposition is realized as the first unit circle
  (for $n=0$ in the definition of the projective limit above). The second component is the
  subgroup $\textbf{Z}_2 \subset\what{Q}_2$ which consists of the elements which have the second
  coordinate in the decomposition equal to $1\in \R/\Z$.
This decomposition gives the coordinatization of the group
$\what{Q}_2$.

Because the group $Q_2$ of dyadic rational numbers can be
represented as the group of all finite on both sides two-sided
sequences of 0's and 1's with the usual binary expansion over $2^n,
n \in \Z$, one can think that the analog of such decomposition for
the group $\what{Q}_2$ is also true. Moreover, we have used
one-sided sequences  with positive indices for parametrizations of
the elements of the subgroup $\textbf{Z}_2$, and that
parametrization agrees with the group structure of 2-adic integers.
Thus, it is tempting to consider the whole group $\what{Q}_2$ of
characters of the group $Q_2$
  as a compact space of all two-sided infinite $\{0,1\}$ sequences:
  $\textbf{X}=\prod_{-\infty}^{+\infty}\{0,1\}=\{0,1\}^\Z$.
 But this is not correct, because there is no needed group structure
 on the space $\textbf{X}$. Nevertheless, it is possible to define a map
  $\pi:\textbf{X}\rightarrow \what{Q}_2$ with the help of the usual dyadic decomposition of the points
  of the unit interval  $(0,1)$ as follows: let $\{x_n\},\, n \in \Z,$
  be a point of $\textbf{X}$; define the pair $(y,\lambda)$, where $y \in \textbf{Z}_2$
  is generated by the sequence
  $\{x_n, n\ge 0\}$ as in Section 1, $$\lambda = \sum_{n=1}^{\infty}  x_{-n}2^{-n}.$$
Denote this map
  by $\pi$:
\begin{equation} \label{def-pi}
\pi:\, \prod_{-\infty}^{+\infty}\{0,1\}\longrightarrow \what{Q}_2,\ \ \ \pi:\,\{x_n\} \mapsto (y,\lambda).
\end{equation}
The map $\pi$ is not an isomorphism of the groups or even
topological spaces but trivially {\em is} an isomorphism (mod 0) of
the measure spaces, where the measure on the space ${\bf X}$ is the
$(1/2,1/2)$ Bernoulli (product) measure, and on the group
$\what{Q}_2,$ it is the Haar measure. So, if we ignore the group
structure of $\what{Q}_2$  and consider it not as solenoid but as a
symbolic space with measure-preserving transformations (odometer,
Morse, etc.), then it is convenient to use the canonical map
$\pi;\,\prod_{-\infty}^{+\infty}\{0,1\}\rightarrow  \what{Q}_2$,
which identifies  only countably many pairs of points. Roughly
speaking, we can consider the 2-solenoid  $\what{Q}_2$ as the space
${\bf X}$ of all two-sided sequences 0's and 1's after some
identifications of elements from the negative (left) side, which
(identification) corresponds to the non-uniqueness of dyadic
decomposition on the left side.

 \subsection{Some transformations and differentiation on the solenoid.}

There is a canonical automorphism $\what S$ on the group
$\what{Q}_2$ : the multiplication by $2$;
 it is conjugate to the automorphism $S^*$ of the group $Q_2$ ---  the multiplication by $1/2$. The
 transformation $\what S$ is a hyperbolic automorphism of the solenoid and in the usual
 coordinatization it is just the Bernoulli 2-shift and {\it natural extension} in Rokhlin
 sense \cite{R} of the one-sided shift $S$  of the space ${\bf Z}_2$.
 which was defined in he section 2.


Now we define the two-sided version of 2-odometer. Let $1$ be the  unit of
the ring $\textbf{Z}_2$ (unity of the multiplicative group).
We extend the odometer $T$ from Section 2, using notation ${\what
T}$, by the same formula $${\what T}x=x+1,$$ where $x$ now is an
element of $\what{Q}_2$. It is useful
to keep in mind that $1$ is a character of the group $Q_2$ which
sends integers $\Z \subset Q_2$ to 1.

 The action of ${\what T}$ does not change the second (left)component in the decomposition
 $\what{Q}_2=Z_2\times \R/\Z$, so it is indeed an extension of $T$. Note that
${\what T}$ is not an ergodic transformation of  $\what{Q}_2$, whereas  $T$ is ergodic on ${\bf Z}_2$.
  We can define also the family of odometer-transformations $$T_0 := \what{T},\ \ T_i:=S^i T_0 S^{-i},\ \ i\in \Z.$$
 It is clear that $T_i$ and $T_j$ commute and its joint action on $\what{Q}_2$ is the action
 of the group  $Q_2$ on $\what{Q}_2$; this is the translation on the corresponding elements, as mentioned above.
   We claim that
\begin{equation} \label{eq-action}
T_i^2 = T_{i+1},\ i\in \Z.
\end{equation}
Indeed, this is immediate for $i=0$ and hence for all $i$.

    Together with the shift $\what S$ the odometers $T_i$ generate a solvable group (wreath product)
     $\Z \rightthreetimes \sum_{\Z} \Z$; the action of this group on the group $\what{Q}_2$
 is continuous and local-transversal in the sense of the paper  \cite{Ver.tran}.

  Define the {\it differentiation} $\what D$ as a transformation of
${\bf X}$ which extends the map $D$ to the space of two-sided sequences.

$$
 {\what D}(\{x_n\}_{-\infty}^{+\infty})=\{(x_n-x_{n+1})\ (\mbox{mod 2})\}
$$

We can, of course, define the differentiation on the solenoid $\what Q_2$ by $\wtil{D} = \pi \circ \what{D} \circ \pi^{-1}$, which is well-defined almost everywhere.
Observe that $\what{D}$ identifies a two-sided sequence
$\til x$ with its ``flip'' $\ov{\til{x}}$, and hence almost everywhere on $\what{Q}_2$ we have
$$
\wtil{D}(y,\lam) = \wtil{D}(z,\gam) \ \Longleftrightarrow\ (y,\lam) = (z,\gam)\ \mbox{or}\ (y ,\lam) = (\ov{z}, 1-\gam).
$$
As we mentioned above, it is difficult to give a precise formula for $\wtil D $ in terms of characters.

  \subsection{Extension of the Morse transformation.}

Now we would like to extend the Morse transformation $M$ from the subgroup ${\bf Z}_2$
to the whole group $\wQ_2$ and space $\textbf{X}$.

We want to have the following properties of $\what M$: it must be a 2-extension of the extended
odometer $\what T$, namely, the relation generalizing Proposition 2.1 must be valid:
\begin{equation} \label{eq-new1}
\what T \circ \what D = \what D \circ \what M,
\end{equation}
and it should be an extension: $$\what{M}|_{\textbf{Z}_2} = M.$$
\begin{theorem}
There is a unique transformation of the space ${\bf X}$
which satisfies the last two equations. Then it defines a measure-preserving transformation $\wtil{M}$ on $\what Q_2$ via
$\wtil{M} = \pi \circ \what{M}\circ \pi^{-1}$, where $\pi$ is defined by (\ref{def-pi}).

\end{theorem}

{\em Proof.}
The uniqueness is clear, and the existence can be shown as follows.
The sequences of the space $\textbf{X}$ can be divided into positive and negative parts:
for $\hat x=(\dots x_{-1},x_0,x_1 \dots)$ denote $x_-=(\dots x_{-2},x_{-1})$ and $x_+=(x_0,x_1 \dots)$. Then we can define
\begin{equation} \label{eq-ext2}
\what{M}(\hat x)\equiv\what{M}((x_-,x_+)) = \left\{\begin{array}{lll}
 (x_-, M(x_+)), & \mbox{if} \ \ \phi(x_+)=0
& \\         (\ov{x_-},M(x_+)), & \mbox{if} \ \ \phi(x)_+=1.
\end{array} \right.
\end{equation}
Here $\phi$ is the cocycle defined in (\ref{def-cocycle}). Verification of (\ref{eq-new1}) is immediate.

On the solenoid ${\what Q}_2$ we get an explicit formula for the Morse transformation:
$$
{\wtil M}(y,\lam) = (My, \lam)\ \ \mbox{if}\ \phi(y) = 0,\ \ \ {\wtil M}(y,\lam) = (Ma, 1 -\lam)\ \ \mbox{if}\ \phi(y) = -1.
$$
(Here $y$ corresponds to $x_+$ in $\bf X$.)
\qed

\medskip

Observe that $\what{M}$ is not ``local'' in the sense that it does change negative coordinates
when the cocycle does not vanish.

The extension $\what{M}$ is continuous on $\{x\in \what{Q}_2:\ x_+\not\in \wtil{\Max}\}$,
which is a set of full (Haar) measure.

Denote $\what{M} = M_0$ and define $M_i=S^i M_0 S^{-i}$; clearly, we have
\begin{equation} \label{eq1}
T_i \circ \what{D} = \what{D} \circ M_i,
\end{equation}
because $\what{D}$ commutes with $S$.

\begin{theorem} The group of transformations, generated by $M_i, i \in \Z$,
is algebraically isomorphic to the group $Q_2$:
\begin{equation} \label{eq-commute}
M_{i+1} = M_i^2,\ \ i\in \Z.
\end{equation}
 Thus, we get a new (Morse) action of
the group $Q_2$ on $\what{Q}_2$. For all $i$, $M_i$ is a 2-point extension of $T_i$.
\end{theorem}

{\em Proof.} We only need to check (\ref{eq-commute}); other statements follow immediately.
Using $S\what{D} = \what{D}S$ and $T_{i+1} = T_i^2$, we obtain that $\what{D} M_{i+1} = \what{D} M_i^2$. It remains to
observe that $M_{i+1}(\hat{x})$ and $M_i^2(\hat{x})$ are cofinal (agree sufficiently far to the right) for all $x\not\in
\wtil{\Max}$. \qed

\medskip

{\it We have defined two canonical measure-preserving actions of the
solvable group $\Z\rightthreetimes Q_2$ on ${\what Q}_2$ --- the
first is generated by the odometer (this is an algebraic action),
and the second which is generated by the Morse action.} Remember
that the Morse action is continuous only almost everywhere.

\medskip

 \textbf{Questions.}

1. Find the cocycle which defines the Morse action as a 2-extension
of the algebraic action analogously to formula (5).

2. Give formula analogous to the formula of Theorem 2.7 which gives
the Morse action on the solenoid as a time change of the algebraic action.

3. How can we characterize both actions of the group
$\Z\rightthreetimes Q_2$ in an intrinsic way?


\begin{thebibliography}{99}
\bibitem{AlSh} J.-P. Allouche and J. Shallit, {\em Automatic sequences.
Theory, applications, generalizations}, Cambridge University Press, Cambridge, 2003.

\bibitem{A} V.I.Arnold. Complexity of finite sequences of zero and
ones and geometry of finite space of functions
{\em Funct.\ Anal.\ Other Math.} {\bf 1} (2006), no.\ 1, 1-15.

\bibitem{B} R.Belinskaya. Partitions of Lebesgue space on trajectories defined by ergodic automorphisms.
{\em Func.\ Anal.\ and its Appl.} {\bf 2} 1968, no.\ 1, 190-199.

\bibitem{Dur} F.~Durand, B.~Host, and C.~Skau,
Substitutional dynamical systems, Bratteli diagrams and dimension   groups,
{\em Ergodic Theory Dynam.\ Syst.}, {\bf 19} (1999), 953--993.

\bibitem{Forrest} A.~Forrest, $K$-groups associated with substitution minimal systems,
{\em Isr.\ J.\ Math.}, {\bf 98} (1997), 101--139.

\bibitem{Kaku} S. Kakutani, Ergodic theory on shift transformations,
{\em Proc.\ 5th Berkeley Symp.\ Math.\ Statist.\ Probab. II}, Berkeley, CA, 1967, pp. 405--414.

\bibitem{Moshe} Y. Moshe, On the subword complexity of Thue-Morse polynomial extractions,
{\em Theoret.\ Comput.\ Sci.} {\bf 389} (2007), no.\ 1-2, 318--329.

\bibitem{Pyth} N. Pytheas-Fogg,
{\em Substitutions in Dynamics, Arithmetics and Combinatorics},
Lecture Notes in Math. {\bf 1794}, Springer-Verlag, 2002.


\bibitem{Quef} M. Queffelec, {\em Substitution Dynamical Systems - spectral analysis}, Springer Lecture Notes in
Mathematics, {\bf 1294} (1987) 1--240.

\bibitem{R} V. A. Rokhlin,
Exact endomorphisms of a Lebesgue space  (Russian),
{\em Izv.\ Akad.\ Nauk SSSR Ser.\ Mat.} {\bf 25} 1961, 499--530.


\bibitem{VerDan} A. M. Vershik, Uniform algebraic approximations of shift and multiplication operators.
{\em Dokl. Akad. Nauk SSSR 259, No.3, 526-529 (1981)}. English translation: {\em Sov. Math. Dokl. 24, 97-100 (1981).}

\bibitem{Ver.tran} A. M. Vershik,
Locally transversal symbolic dynamics, (Russian)  {\em Algebra i Analiz} {\bf 6}  (1994),  no.\ 3, 94--106;
English translation:
{\em St.\ Petersburg Math.\ J.} {\bf 6}  (1995),  no.\ 3, 529--540.

\bibitem{VerLiv} A. M. Vershik, A.N. Livshits, Adic models of ergodic
transformations, spectral theory, and related topics,
{\em Adv.\ in Soviet Math. AMS Transl.} {\bf 9} (1992), 185--204.

\end{thebibliography}
\end{document}